\documentclass[a4paper]{amsart}

\usepackage{amssymb,latexsym}
\usepackage{amsmath}
\usepackage{amsfonts}
\usepackage{amsthm}
\usepackage{amssymb}

\theoremstyle{plain}
\newtheorem{theorem}{Theorem}

\newtheorem{proposition}[theorem]{Proposition}

\theoremstyle{definition}

\title{Diophantine quadruples with the properties $D(n_1)$ and $D(n_2)$}
\begin{document}

\date{}


\author[A. Dujella]{Andrej Dujella}
\address{
Department of Mathematics\\
Faculty of Science\\
University of Zagreb\\
Bijeni{\v c}ka cesta 30, 10000 Zagreb, Croatia
}
\email[A. Dujella]{duje@math.hr}

\author[V. Petri\v{c}evi\'c]{Vinko Petri\v{c}evi\'c}
\address{
Department of Mathematics\\
Faculty of Science\\
University of Zagreb\\
Bijeni{\v c}ka cesta 30, 10000 Zagreb, Croatia
}
\email[V. Petri\v{c}evi\'c]{vpetrice@math.hr}

\begin{abstract}
For a nonzero integer $n$, a set of $m$ distinct nonzero integers $\{a_1,a_2,\ldots,a_m\}$
such that  $a_ia_j+n$ is a perfect square for all $1\leq i<j\leq m$, is called a $D(n)$-$m$-tuple.
In this paper, we show that there infinitely many essentially different quadruples
which are simultaneously $D(n_1)$-quadruples and $D(n_2)$-quadruples with $n_1\neq n_2$.
\end{abstract}


\maketitle

\section{Introduction}

For a nonzero integer $n$, a set of distinct nonzero integers $\{a_1,a_2,\ldots,a_m\}$ such that
$a_ia_j+n$ is a perfect square for all $1\leq i<j\leq m$, is called a  Diophantine $m$-tuple
with the property $D(n)$ or $D(n)$-$m$-tuple.
The $D(1)$-$m$-tuples are called simply Diophantine $m$-tuples,
and have been studied since the ancient times.
Diophantus of Alexandria found a set of four rationals
$\left \{\frac{1}{16},\frac{33}{16},\frac{17}{4},\frac{105}{16}\right \}$ with
the property that the product of any two of its distinct elements is a square of a rational number.
By multiplying elements of this set by $16$ we obtain the $D(256)$-quadruple
$\{1,33,68,105\}$. Fermat found the first $D(1)$-quadruple, it was the set $\{1,3,8,120\}$.
In 1969, Baker and Davenport \cite{BD69}, using linear forms in logarithms of algebraic numbers and
the reduction method introduced in that paper,
showed that the set $\{1,3,8\}$ can be extended to a Diophantine quintuple only by adding $120$ to the set.
In 2004, Dujella~\cite{d05} showed that there are no Diophantine sextuples and that there are at most finitely many Diophantine quintuples.  Recently, He, Togb\'e and Ziegler proved that there are no Diophantine quintuples~\cite{HTZ16+}.
(See also \cite{BTF} for an analogous result concerning the conjecture of nonexistence of $D(4)$-quintuples.)
On the other hand, it was known already to Euler that there are infinitely many rational Diophantine quintuples. In particular, the Fermat's set $\{1,3,8,120\}$ can be extended to a rational Diophantine quintuple by adding
$777480 / 8288641$ to the set. Recently, Stoll \cite{Stoll} proved that the extension of Fermat's set to a rational
Diophantine quintuple is unique. The first example of a rational Diophantine sextuple, the
set $\{11/192, 35/192, 155/27, 512/27, 1235/48, 180873/16\}$, was found by Gibbs \cite{Gibbs1},
while Dujella, Kazalicki, Miki\'c and Szikszai \cite{DKMS} recently proved that there are infinitely
many rational Diophantine sextuples (see also \cite{Duje-Matija}). It is not known whether there
exists any rational Diophantine septuple. For an overview of results on $D(1)$-$m$-tuples and its generalizations see \cite{Duje-Notices}.

Let us mention some results concerning $D(n)$-sets with $n\neq 1$.
It is easy to show that there are no $D(n)$-quadruples if $n \equiv 2 \pmod{4}$ (see e.g. \cite{Brown}).
On the other hand, it is known that if $n\not\equiv 2 \pmod{4}$ and $n \not\in \{-4,-3,-1, 3, 5, 8, 12, 20\}$, then there exists at least one $D(n)$-quadruple \cite{Duje-acta}.
It is believed that the size of sets with the property $D(n)$ is bounded by an absolute constant
(independent on $n$). It is known that the size of sets with the property $D(n)$
is $\leq 31$ for $|n|\leq 400$; $< 15.476 \log|n|$ for $|n| > 400$, and $ < 3 \cdot 2^{168}$ for $n$ prime (see \cite{Duje-size1,Duje-size2,DL-IMRN} and also \cite{murty}).

In \cite{kk01}, A.~Kihel and O.~Kihel asked if there are Diophantine triples  $\{a, b, c\}$
which are $D(n)$-triples for several distinct $n$'s. They conjectured that  there are no
Diophantine triples which are also $D(n)$-triples for some $n\neq 1$.
However, the conjecture is not true, since, for example, $\{8, 21, 55\}$ is a $D(1)$ and $D(4321)$-triple
(as noted in the MathSciNet review of \cite{kk01}),
while $\{1, 8, 120\}$ is a $D(1)$ and $D(721)$-triple, as observed by Zhang and Grossman~\cite{zg15}.
In \cite{ADKT}, several infinite families of Diophantine triples were presented
which are also  $D(n)$-sets for two additional $n$'s.
Furthermore, there are examples of Diophantine triples which are $D(n)$-sets for three additional $n$'s.
For example, the set $\{4, 12, 420\}$ is a $D(n)$-quadruple for $n=1, 436, 3796, 40756$ 
(see also \cite{ADKT-rims}).

\medskip

If we omit the condition that one of the $n$'s is equal to $1$,
then the size of a set $N$ for which there exists a triple $\{a, b, c\}$ of nonzero integers
which is a $D(n)$-set for all $n\in N$ can be arbitrarily large. Indeed,
take any triple $\{a, b, c\}$ such that the elliptic curve
$$ E \,:\, \quad y^2 = (x+ab)(x+ac)(x+bc) $$
has positive rank over $\mathbb{Q}$. Then there are infinitely many rational points
on $E(\mathbb{Q})$. For an arbitrary large positive integer $M$ we may choose $M$ distinct rational points $R_1,\ldots,R_M\in 2E(\mathbb{Q})$, so that we have
\[
x(R_i)+bc=\square,\quad  x(R_i)+ca=\square, \quad x(R_i)+ab =\square,
\]
where $\square$ stands for a square of a rational number (see e.g. \cite[4.1, p.~37]{hu}).
We choose $z\in \mathbb{Z}\setminus\{0\}$ such that $z^2x(R_i)\in \mathbb{Z}$ for all $i=1, 2,\ldots, m$.
Then the triple $\{az, bz, cz\}$ is a $D(n)$-triples for
$n=x(R_i)z^2$ for all $i=1,2,\ldots,m$ (see \cite[Section 4]{ADKT} for the details).

On the other hand, assuming Lang's conjecture on varieties of general type,
for a given quadruple $\{a, b, c,d\}$ of distinct integers,
the size of the set $N$ of integers $n$ for which $\{a, b, c,d\}$ is a $D(n)$-quadruple
is bounded by an absolute constant. Indeed, let $ab+n=x^2$. By multiplying remaining five conditions,
we get the hyperelliptic curve
$$ y^2 = (x^2+ac-ab)(x^2+bc-ab)(x^2+ad-ab)(x^2+bd-ab)(x^2+cd-ab), $$
which has genus $4$ unless it has two equal roots. Assume e.g. that $ad=bc$.
Then we get the hyperelliptic curve
$$ y^2 = (x^2-ab+ac)(x^2+bc-ab)(ax^2-a^2b+b^2c)(ax^2-a^2b+bc^2) $$
with distinct roots (unless $b=-a$ or $c=-a$) and, hence, with genus equal to $3$.
Finally, if e.g. $c=-a$ and $d=-b$, we get the hyperelliptic curve
$$ y^2 = (x^2-ab-a^2)(x^2-2ab)(x^2-ab-b^2) $$
with distinct roots and with genus equal to $2$.
Assuming the above mentioned Lang' conjecture, Caporaso, Harris and Mazur \cite{C-H-M}
proved that for $g\ge 2$ the number  $B(g,\mathbb{K})=\max_{C} |C(\mathbb{K})|$ is finite,
where $C$ runs over all curves of genus $g$ over a number field $\mathbb{K}$.
Therefore, we get that, under Lang's conjecture, $|N| \leq \max(B(2,\mathbb{Q}),B(3,\mathbb{Q}),B(4,\mathbb{Q}))$.

\medskip

Thus, it seems natural to ask is there any set of four distinct nonzero integers
which is a $D(n_i)$-quadruple for two distinct (nonzero) integers $n_1$ and $n_2$.
However, it seems that this question has not been studies yet and that there are no examples
of such quadruples in the literature. In Section \ref{sec:numex} we will present results
of our computer search for such quadruples. Motivated by certain regularities in found examples,
we will show in Section \ref{sec:17} that there are infinitely many such examples.
If $\{a,b,c,d\}$ is $D(n_1)$ and $D(n_2)$-quadruple and $u$ is a nonzero rational such that
$au,bu,cu,du,n_1u^2$ and $n_2u^2$ are integers, then $\{au,bu,cu,du\}$ is a $D(n_1u^2)$
and $D(n_2u^2)$-quadruples. We will say that these two quadruples are equivalent,
and list only one representative of each found class of quadruples.

Our main result is

\begin{theorem} \label{tm:1}
There are infinitely many nonequivalent sets of four distinct nonzero integers
$\{a,b,c,d\}$ with the property that
there exist two distinct nonzero integers such that $\{a,b,c,d\}$ a $D(n_1)$-quadruple and
a $D(n_2)$-quadruple.
\end{theorem}

\section{Numerical examples} \label{sec:numex}

We started with computational search for $D(n)$-quadruples, where $-500\,000\le n\le 500\,000$.
For a fixed nonzero integer $n$, by observing divisors of integers of the form $m^2-n$,
it is not hard to get some $D(n)$-quadruples (we were searching in the range $m\le 333\,333$).

We have implemented the algorithm in {\tt C++}.
For a fixed $n$, we construct a graph, connecting the numbers $k$ and $l$ with an edge
provided they satisfy $ k\cdot l=m^2-n$.
The graph can be represented using standard containers
(for example {\tt map<long, set<long> > g}; so for $k<l$, $k$ and $l$ are connected if set {\tt g[l]} contains $k$).
We also connect $k$ and $l$ with $k+l+2m$, since $k(k+l+2m)+n=(k+m)^2$ and
$l(k+l+2m)+n=(l+m)^2$ (a $D(n)$-triple of the form $\{k,l,k+l\pm 2\sqrt{kl+n}\}$ is called
\emph{regular}).

But we actually used container {\tt unordered\_map<long, vector<long> >},
which is somewhat faster and takes less memory. For $m=1,\dots, 333\,333$,
it usually takes about 10--12 seconds (on one core of 3.6GHz)
to build such a graph, and it usually takes about 500MB of memory
(but graph density depends on $n$). Then we search for a 4-clique in graph (e.g. $D(n)$-quadruple).
We do this by sorting each {\tt vector}, and using binary search.
So for finding all 4-cliques it takes about a second, and for the most of $n$'s
we get several hundreds of quadruples.

Then we searched for $n_2$ using M. Stoll's program {\tt ratpoints} (see \cite{ratpoints}).
For a quadruple $\{a,b,c,d\}$, the search for an integer point on the hyperelliptic curve $y^2=(ab+x)(ac+x)(ad+x)(bc+x)(bd+x)(cd+x)$ with $x=n_2\le10^8$ takes about $0.02$ seconds.
Here is summarize results of our search:

\begin{center}
\begin{tabular}{r|l@{\quad\quad}r|l}
	$\{a,b,c,d\}$&$\{n_1,n_2\}$ 		&$\{a,b,c,d\}$&$\{n_1,n_2\}$\\
	\hline
-1701, -901, 224, 243 &413424, 463968  &-1, 7, 22532, 23407 *&30632, 214376\\
-189, -133, 27, 32 *&6192, 8352	       &15, 380, 5735, 634880 &361536, 7123200\\
-176, -169, 169, 176 &31265, 36305	   &15, 720, 9135, 40656 &17424, 13708816\\
-52, 135, 351, 575 &37296, 67536	   &27, 115, 160, 1755 &-2016, 37296\\
-27, 28, 189, 493 *&13752, 61272       &28, 6348, 18750, 88872 &330625, 38101225\\	
-27, 189, 4189, 6364 *&194328, 1325304 &45, 276, 8820, 18228 &112896, 2966656\\
-15, 1140, 2057, 15609 &234256, 989296 &51, 192, 315, 2331 &-6656, 1080144\\
-11, 28, 385, 540 &11124, 34164		    &69, 300, 949, 2925 &63400, 417544\\
-4, 209, 5129, 49049 &252840, 6062280   &70, 430, 2178, 18634 &-20691, 1678149\\	
-3, 21, 1152, 1517 *&5392, 37312		&125, 2709, 2816, 5621 &-273600, 1443600\\
-3, 21, 2597, 3132 *&11512, 80152		&169, 448, 8640, 11305 &97344, 28482624\\
-1, 7, 64, 119 *&128, 848				&175, 231, 300, 396 &-16400, -40400\\	
-1, 7, 4484, 4879 *& 6248, 43688		&234, 322, 406, 1222 &-10323, -69723\\
	\end{tabular}
	\end{center}

We indicate by * quadruples which contain two elements $a$ and $b$ such that $a/b=-1/7$.
These quadruples will play crucial role in the proof of Theorem \ref{tm:1} in the next section.

\section{Quadruples containing the pair $\{-1,7\}$} \label{sec:17}

Motivated by the examples indicated by * in the previous section,
we will show that there infinitely many quadruples of the form
$\{a,b,c,d\}$, where $a/b=-1/7$ that are $D(n)$-quadruples for two distinct (nonzero) $n$'s.
Then we will show that in fact we may take $a=-1$ and $b=7$ and get the the same conclusion.

We use regular triples mentioned in the previous section.
Namely, if $AB+n=R^2$, then $\{A,B,A+B+2R\}$ and
$\{A,B,A+B-2R\}$ are $D(n)$-triples.
Let $cd+n_1=r^2$ and $cd+n_2=s^2$.
If $c+d-2r=7$ and $c+d-2s=-1$, then $\{7,c,d\}$ is a $D(n_1)$-triple and $\{-1,c,d\}$ is a
$D(n_2)$-triple. We have to satisfy the remaining six conditions from the definition of
$D(n_i)$-quadruples.

We search for a solution in the from
$n_2=kn_1-l$ with (rational) constants $k$ and $l$.
We have $n_1= -cd+r^2$ and $c = -d+2r+7$. From $c+d-2s=-1$, we get $r = -4+s$.
By inserting this in $cd+n_2=s^2$ and solving this equation for $d$, we get that
$(28k^2-24k-4)s-63k^2+62k+4kl-4l+1$ is a perfect square, say $t^2$. Thus we obtain
\begin{align*}
s &=\frac{4l+63k^2-62k-4kl-1+t^2}{4(7k^2-6k-1)}, \\
d &=\frac{-4kl+4l+1-50k+49k^2-2t-14tk+t^2}{4(7k^2-6k-1)}.
\end{align*}
Consider now the condition that $bd+n_2$ is a perfect square.
We obtain a quadratic function in $t$ with the discriminant $-128(k-1)^2(112k^2-112k-15kl+7l)$.
From $112k^2-112k-15kl+7l=0$ we get
$$ l = \frac{112k(k-1)}{15k-7}. $$
There are four remaining conditions: $ab+n_2$, $ab+n_1$, $ac+n_1$ and $ad+n_1$ are perfect squares.
All four conditions lead to quadratic functions in $t$. The corresponding discriminants are
\begin{gather*}
56k(15k-7)(7k+1)(13k-7)(k-1)^2, \\
-56(15k-7)(7k+1)(k-7)(k-1)^2, \\
2048(7k+1)(k-7)(k-1)^3(15k-7), \\
2048(7k+1)(k-7)(k-1)^3(15k-7)
\end{gather*}
(last two discriminants are identical).
Hence, by taking $k=7$, we can satisfy last three conditions simultaneously.
Only one condition remains, $ab+n_2$ is a perfect square, and it is equivalent to $-6+\frac{7}{900}t^2$
being perfect square. From
$$ -6+\frac{7}{900}t^2= (1+ \frac{u(t-30)}{30})^2, $$
we obtain
\begin{align*}
c &= \frac{(2u^2-3u+7)(2u^2-u-7)}{(u^2-7)^2},\\
d &= -\frac{(u^2-3u+14)(u^2+u-14)}{(u^2-7)^2}, \\
n_1 &= \frac{4(2u^4-u^3-20u^2-7u+98)}{(u^2-7)^2}, \\
n_2 &= \frac{4(2u^2-7u+14)(u^2+7)}{(u^2-7)^2}.
\end{align*}
For $u\not \in \{0,1,2,-7/5,-5,7/2,7,4,7/3,-7,-2,3,-7/2,7/4\}$ the elements of the set
$\{-1,7,c,d\}$ are distinct rationals. By taking $u=v/w$ and getting rid of denominators,
we obtain the following result.

\begin{proposition} \label{prop:1}
Let $v$ and $w$ be coprime integers and
$$ v/w \not\in \{0,1,2,-2,3,4,-5,7,-7,7/2,-7/2,7/3,7/4,-7/5\}. $$
Then the set
\begin{equation}\label{eq:vw}
\begin{gathered}
 \{ -(-v^2+7w^2)^2, 7(-v^2+7w^2)^2, -(-2v^2+vw+7w^2)(2v^2-3vw+7w^2), \\
 (v^2-3vw+14w^2)(-v^2-vw+14w^2) \}
\end{gathered}
\end{equation}
is a $D(n_1)$-quadruple and a $D(n_2)$-quadruple for
\begin{align*}
n_1&=4(-v^2+7w^2)^2(2v^4-v^3w-20v^2w^2-7vw^3+98w^4), \\
n_2&=4(-v^2+7w^2)^2(2v^2-7vw+14w^2)(v^2+7w^2).
\end{align*}
\end{proposition}

We have obtained infinitely many quadruples with the required property satisfying $a/b=-1/7$
(in other word, infinitely many rational quadruples with $a=-1$, $b=7$).
Now we will show that there are infinitely many integer quadruples with $a=-1$, $b=7$.
Indeed, let $v$ and $w$ be a solution of the Pellian equation
\begin{equation} \label{eq:pell}
v^2 - 7w^2 = 2.
\end{equation}
The equation (\ref{eq:pell}) has infinitely many integer solutions given by
$$ v_0=3, \quad v_1=3, \quad v_{i+2}=16v_{i+1}-v_i, $$
$$ w_0=-1, \quad w_1=1, \quad w_{i+1}=16w_{i+1}-w_i. $$
By inserting $v=v_i$, $w=\pm w_i$ in (\ref{eq:vw}), and dividing elements of the quadruple
by the common factor $4$, we obtain quadruples of the form $\{-1,7,c,d\}$ which
are $D(n)$-quadruples for two distinct $n$'s. Here are few smallest examples:

\begin{center}
\begin{tabular}{r|l}
	$\{a,b,c,d\}$&$\{n_1,n_2\}$ \\		
	\hline
-1, 7, 119, 64 & 128, 848 \\
-1, 7, 1191959, 1185664 & 1585088, 11095568 \\
-1, 7, 5840864, 5826919 & 7778528, 54449648 \\
-1, 7, 76695715424, 76694116519 & 102259887968, 715819215728 \\
-1, 7, 376369378007, 376365836032 & 501823476032, 3512764332176
	\end{tabular}
	\end{center}

\section{The case $n_1=0$}

In the definition of $D(n)$-$m$-tuples, the case of $n=0$ is usually excluded,
although certainly the definition make sense in this case also.
The reason for excluding $n=0$ is in very different behavior of $D(0)$-tuples compared
with $D(n)$-tuples for $n\neq 0$. While for a fixed $n\neq 0$ the size of sets with the property $D(n)$
is bounded, sets with the property $D(0)$ can be arbitrarily large, just take any subset of the set of squares
$\{1,4,9,16,\ldots\}$. However, in the context of finding quadruples which are $D(n_1)$ and $D(n_2)$-quadruples
for $n_1\neq n_2$, it seems to be natural to consider also the case $n_1=0$.
We might expect that in this case it could be easier to find such quadruples,
but it seems that there is not straightforward way to see why there should be infinitely many of them.

A simple search for $D(n)$-quadruples which elements are perfect squares gives many such examples.
Here we list some of them:

\begin{center}
\begin{tabular}{r|l@{\qquad\quad}r|l}
	$\{a,b,c,d\}$&$\{n_1,n_2\}$ 		&$\{a,b,c,d\}$&$\{n_1,n_2\}$\\
	\hline
1, 4, 169, 1024 &0, 6720	        &196, 625, 1024, 3969 &0, 705600\\
1, 36, 529, 1024 &0, 60480	        &324, 841, 1369, 4096 &0, -262080\\
25, 64, 961, 2025 &0, 188496	    &1, 324, 2209, 4096 &0, 887040\\	
100, 625, 1024, 2025 &0, 360000		&36, 729, 2500, 4096 &0, 518400 \\
64, 169, 441, 2401 &0, 1164240	    &256, 729, 2401, 5625 &0, 1587600 \\
961, 1849, 2704, 2916 &0, -1774080  &121, 169, 2704, 5625 &0, 1436400 \\
32, 98, 1152, 3528 &0, 257985       &1681, 4096, 5625, 5929 &0, -6879600
	\end{tabular}
	\end{center}

Our starting point in constructing infinitely many quadruples $\{a,b,c,d\}$ which are
$D(0)$-quadruples and also $D(n_2)$-quadruples for $n_2\neq 0$, is the following simple fact
(see \cite[Theorem 1]{Duje-graz} and \cite[Section 5]{Duje-acta}). The set
$$ \{ a, ak^2-2k-2, a(k+1)^2-2k, a(2k+1)^2-8k-4 \} $$
is a $D(2a(2k+1)+1)$-quadruple provided all its elements are distinct nonzero integers.
Thus, we take $b=ak^2-2k-2$, $c=a(k+1)^2-2k$, $d=a(2k+1)^2-8k-4$, $n=2a(2k+1)+1$,
and we want to find integers $a$ and $k$ such that $\{a,b,c,d\}$ is also a $D(0)$-quadruple,
i.e. such that $ab$, $ac$ and $ad$ are perfect squares.
By putting $ab=(ak+r)^2$ we get
$$ k = -\frac{2a+r^2}{2(a(1+r))}. $$
Then we put $ac=(2ra+s)^2$ and we get
$$ a = -\frac{-4r^2-4r^3-r^4+s^2}{4(r^3-2-2r+rs)}. $$
The final condition that $ad$ is a perfect square, now becomes
\begin{equation}\label{eq:rsn0}
\begin{aligned}
&(r^2-2r+1)s^4+(8r^4+8r^3-16r)s^3 \\
&\mbox{}+ (22r^6+68r^5+54r^4-40r^3-24r^2+32r+32)s^2 \\
&\mbox{}+ (-192r^3+24r^8-448r^4+88r^7-336r^5)s \\
&\mbox{}+9r^{10}+30r^9-39r^8-248r^7-200r^6+352r^5+752r^4+512r^3+128r^2 =\Box.
\end{aligned}
\end{equation}
Since $r^2-2r+1$ is a square, this quartic curve in $s$ has rational point at infinity,
so it can be in standard way transformed into an elliptic curve over $\mathbb{Q}(r)$:
\begin{equation}\label{eq:xyn0}
\begin{aligned}
y^2 &= x^3+(4r^6+56r^5+84r^4+80r^3+48r^2-64r-64)x^2\\
&\mbox{}+ (-1024r^9-2048r^8+1024r^7+5120r^6+3072r^5-3072r^4\\
&\mbox{}-5120r^3-1024r^2+2048r+1024)x.
\end{aligned}
\end{equation}
The curve (\ref{eq:xyn0}) has a point $[0,0]$ of order $2$ and two independent points of infinite order:
\begin{align*}
P_1 &=[-6r^6-4r^5+74r^4+168r^3+88r^2-32r-32,\\
& \hspace*{0.5cm}-256r^7-1792r^6-4352r^5-4352r^4-1024r^3+1024r^2+512r], \\
P_2 &= [-4r^5+10r^4+8r^3+24r^2-6r^6-32r, \\
& \hspace*{0.5cm}-320r^3+128r+448r^6-512r^4-64r^2+128r^7+192r^5].
\end{align*}
If fact, by using the algorithm of Gusi\'c and Tadi\'c from \cite{GT2}
(see also \cite{GT1,Stoll} for other variants of the algorithm), we can check that
the rank of (\ref{eq:xyn0}) over $\mathbb{Q}(r)$ is equal to $2$ and that $P_1$ and $P_2$
are its free generators. Indeed, the specialization $r=13$ satisfies the assumptions of
\cite[Theorem 1.3]{GT2}.

Hence, there are infinitely many rational points on curves (\ref{eq:xyn0}) and (\ref{eq:rsn0}),
and thus infinitely many quadruples with the required property.
We present an explicit formula.
By taking the point $P_2-P_1$ on  (\ref{eq:xyn0}) and we get
$$ s=-\frac{r(3r^3+9r^2+7r+2)}{r^2+r-1}, $$
and (after multiplying with the common denominator) the quadruple
\begin{equation} \label{eq:rn0}
 \{4r^4(r+2)^2, (r^3-4r+1)^2, (r^3+4r^2-1)^2, 4(2r-1)^2 \}
\end{equation}
which is a $D(0)$-quadruple and a
$D(16r^{10}+96r^9+112r^8-192r^7-256r^6+192r^5+112r^4-96r^3+16r^2)$-quadruple.
By taking $r$ to be an integer in (\ref{eq:rn0}) we obtain the following result

\begin{proposition} \label{prop2}
There are infinitely many nonequivalent sets of four distinct nonzero integers
$\{a,b,c,d\}$ with the property that $a,b,c,d$ are perfect squares (so that $\{a,b,c,d\}$ is
a $D(0)$-quadruple) and there exist $n_2\neq 0$ such that $\{a,b,c,d\}$ a $D(n_2)$-quadruple.
\end{proposition}

Let us mention that in \cite{DGPT,MacLeod} sets which all elements are squares
appeared in similar context (construction of (strong)
Eulerian $m$-tuples, which are shifted $D(-1)$-$m$-tuples). Other connections of
(rational) Diophantine $m$-tuples and elliptic curves can be found in
\cite{ADP,Duje-Glasnik,DKMS,Duje-Peral-LMS,Duje-Peral-RACSAM}.

\bigskip

{\bf Acknowledgements.}
The authors were supported by the Croatian Science Foundation under the project no.~IP-2018-01-1313.
The authors acknowledge support from the QuantiXLie Center of Excellence, a project
co-financed by the Croatian Government and European Union through the
European Regional Development Fund - the Competitiveness and Cohesion
Operational Programme (Grant KK.01.1.1.01.0004).
The authors acknowledge the usage of the supercomputing resources
of Division of Theoretical Physics at Ru\dj{}er Bo\v{s}kovi\'c Institute.

\end{document}